\documentstyle[11pt]{article}
\begin{document}
\hoffset-2cm \oddsidemargin 1.0in
\title{{\large \bf Test Map and Discreteness Criteria for Subgroups in PU(1,$n$;C)}
\thanks{The Project-sponsored by SRF for ROCS, SEM and NSFC(No.10771200)}
}
\author{\normalsize Chang-Jun Li and Xiao-Yan Zhang}
\date{}
\maketitle \baselineskip 18pt
\begin{minipage}{130mm}
\begin{abstract}
{In this paper, we study the discreteness for non-elementary
subgroup $G$ in $PU(1,n;C)$, under the assumption that $G$
satisfies Condition $A$. Mainly, we present that one can use a
test map, which need not to be in $G$, to examine the discreteness
of G, and also show that $G$ is discrete, if every
two-loxodromic-generator subgroup of $G$ is discrete. }
\end{abstract}
\end{minipage}\\

\vskip 2pt
 \vskip 2mm {\bf 1. Introduction}\vskip
 2mm

The discreteness of M\"{o}bius groups is a fundamental problem,
which have been investigated by a number of authors. In 1976,
J{\o}rgensen [13] proved a necessary condition for a
non-elementary two generator subgroup of $SL(2,C)$ to be discrete,
which is called J{\o}rgensen$^{,}$s inequality. By using this
inequality, J{\o}rgensen established the following famous
result[14]:

 \textbf{
Theorem 1.1.} A non-elementary subgroup G of $SL(2,C)$ is discrete
if and only if all its two-generator subgroups are discrete.

This important result has become standard in literature and it
indicates that the discreteness of a non-elementary M\"{o}bius
group depends on the information of all its rank two subgroups.
Furthermore, Gilman [15] and Isochenko [20] strengthened the above
theorem, and showed that $G$ is discrete if every subgroup
generated by two loxodromic elements is discrete. There are many
further discussions about discreteness criteria in this direction.
For more details, see the references[9,21,22,23].\

\par In [1,7,16,18], the authors have discussed the generalization of
Theorem 1.1 to higher dimensional hyperbolic space. Moreover, Fang
and Nai [7] also obtained the following result:

\textbf{Theorem 1.2.} Let a non-elementary subgroup $G$ of
$SL(2,\Gamma_{n})$ satisfy condition $A$. Then $G$ is discrete if
and only if two arbitrary loxodromic elements $f$ and $g$ in $G$
the group $\langle f,g \rangle$ is discrete.

In 2004, Chen min [19] showed that one can even use a fixed
M\"{o}bius transformations as a test map to test the discreteness
of a given M\"{o}bius group. More precisely,

\textbf{Theorem 1.3.} Let $G$ be an $n$-dimensional subgroup of
Isom($H^{n}$), and $f$ be a non-trivial M\"{o}bius transformation.
If for each $g\in G$, the group $\langle f,g \rangle$ is discrete,
then $G$ is discrete.

The result suggests that the discreteness is not a totally
interior affair of the involved group, and this provides a new
point of view to the discreteness problem.

In complex hyperbolic space, Kamiya [17] established a similar
version of theorem 1.1 for finitely generated subgroups of
$PU(1,n;C)$ as follows:

\textbf{Theorem 1.4.} Let G be a non-elementary finitely generated
subgroups of PU(1,$n$;C), then G is discrete if and only if
$\langle f,g\rangle$ is discrete for any $f$ and $g$ in G.

In 2001, Dai B,Fang and Nai [6] proved that:

\textbf{Theorem 1.5.} Let G be a non-elementary  subgroup of
PU(1,$n$;C) with condition A, then G is discrete if and only if
$\langle f,g\rangle$ is discrete for any $f$ and $g$ in G.

Here, G is said to satisfy \textbf{condition A }if it has no
sequence \{$g_{i}$\} of distinct elements of finite order such
that Card(fix($g_{i}$))=$\infty$ and $g_{i}\rightarrow I$ as
$i\rightarrow \infty$, where $$fix(g_{i})=\{x\in
\partial H^{n}_{C}:g_{i}(x)=x\}.$$

In this paper, we continue to discuss the discreteness criteria
for non-elementary subgroup $G$ in $PU(1,n;C)$, and we will
acquire three conclusions under the assumption that $G$ satisfies
Condition $A$. The first result is similar to Theorem 1.3, which
primarily consider to use a parabolic or loxodromic element as a
test map to examine the discreteness of $G$, but whether one can
use a elliptic element remains a open problem. The next result is
followed from the idea of Theorem 1.2, and it shows that $G$ is
discrete, if each two-loxodromic-generator subgroup of $G$ is
discrete. And the third conclusion strengthened the second result,
for details, see the section 3.

 \vskip 2pt
 \vskip 2mm {\bf 2. Notations and Preliminary Results}\vskip
 2mm

Throughout this paper, we will adopt the same notations and
definitions as in [4,10,12]. Now we start by giving some general
facts about $PU(1,n;C)$.

Let C be the field of complex numbers, $V=V^{1,n}(C)(n\geq1)$
denote the vector space $C^{n+1}$, together with the unitary
structure defined by the Hermitian form
$$\Phi(z^{*},w^{*})=-\overline{z}^{*}_{0}w^{*}_{0}+\sum_{j=1}^{n} -\overline{z}^{*}_{j}w^{*}_{j}$$
for $z^{*}=(z^{*}_{0},z^{*}_{1},...,z^{*}_{n}),
w^{*}=(w^{*}_{0},w^{*}_{1},...,w^{*}_{n})\in V.$

An automorphism $g$ of V, that is a linear bijection such that
$$\Phi(z^{*},w^{*})=\Phi(g(z^{*}),g(w^{*}))$$
for $z^{*},w^{*}\in V$, will be called a unitary transformation.
We denote the group consisting of all unitary transformation by
$U(1,n;C)$. Let
$$V_{0}=\{z^{*}\in V: \Phi(z^{*},z^{*})=0\},\;\;\;\;V_{-}=\{z^{*}\in V:\Phi(z^{*},z^{*})<0\}.$$
Set$$PU(1,n;C)=U(1,n;C)/(center).$$ It is obvious that $V_{0}$ and
$V_{-}$ are invariant under $U(1,n;C)$. Set
$$V^{s}=V_{-}\bigcup V_{0}\setminus \{ 0 \}.$$ Let $P:V^{s}\rightarrow
P(V^{s})$ be the projection map defined by
$$P(z_{0}^{*},z_{1}^{*},...,z_{n}^{*})=(z_{1}^{*}z_{0}^{-1},...,z_{n}^{*}z_{0}^{-1}).$$
We denote $P(0,1,...,0)$ by $\infty$. We may identify $P(V_{-})$
with the Siegel domain
$$H^{n}=\{w=(w_{1},w_{2},...,w_{n})\in C^{n}: \; Re(w_{1})>\frac{1}{2}\sum_{j=2}^{n}|w_{j}|^{2}\}.$$

An element of $PU(1,n;C)$ acts on $H^{n}_{C}$ and its boundary
$\partial H^{n}_{C}.$ Denote $H^{n}_{C}\bigcup \partial H^{n}_{C}$
by $\overline{H^{n}_{C}}.$
As in [4,12], a non-trivial element $g$ in $PU(1,n;C)$ is called\\
(1) elliptic if it has a fixed point in $H^{n}_{C};$\\
(2) parabolic if it has exactly one fixed point and the point lies
on
$\partial H^{n}_{C}$;\\
(3) loxodromic if it has exactly two fixed points and the points
lie on $\partial H^{n}_{C}.$

For a subgroup $G\subset PU(1,n;C)$, the limit set $L(G)$ of $G$
is defined as $$L(G)=\overline{G(p)}\bigcap \partial
H^{n}_{C}(p\in H^{n}_{C}).$$ The fixed point sets of $f\in G$ and
of $G$ are
$$fix(f)=\{x\in \overline{H^{n}_{C}}:f(x)=x\},\;\;fix(G)=\bigcap_{f\in G} fix(f).$$

\textbf{Definition 2.1}[12]. A subgroup $G\subset PU(1,n;C)$ is
said to be non-elementary, if $G$ contains two non-elliptic
elements of infinite order with distinct fixed points, or $G$ is
said to be elementary.

 \textbf{Definition 2.2}[12]. $G_{L}=\{g\in
G: g(x)=x, for \;any\;\; x\in L(G)\}.$

 \textbf{Definition 2.3}[12]. A subgroup
$G\subset PU(1,n;C)$ is said to be bounded torsion if there exists
an integer number $m$ such that for each $g\in G$ has $ord(g)\leq
m$ or $ord(g)=\infty.$

\textbf{Definition 2.4}[12]. Let $X$ be subgroup of the vector
space $V$. The span of $X$ denoted as $\langle X \rangle$ is the
smallest $C$-subspace containing $X$. If $X$ is a subset of
$H^{n}_{C}$, the span $\langle X \rangle$ is defined by $\langle X
\rangle=P(\langle P^{-1}(X)\rangle)\bigcap V_{-}$.

 \textbf{Lemma
2.5} (Lemma2.1 of [3]). Suppose that $f$ and $g\in PU(1,n;C)$
generate a discrete and
non-elementary group. Then\\
i) if $f$ is parabolic or loxodromic, we have
$$max\{N(f),N[f,g])\}\geq 2-\sqrt{3}$$
where $[f,g]=fgf^{-1}g^{-1}$ is the commutator of $f$ and $g$,
$N(f)=\| f-I \|$.\\
ii) if $ f $ is elliptic, we have
$$max \{N(f),N([f,g ^{i }]):i= 1,2,...,n+1  \} \geq
2-\sqrt{3}.$$

\textbf{Lemma 2.6.} Let $G $ be a non-elementary subgroup of
$PU(1,n;C)$ and let $O_{1} $ and $O_{2} $ be disjoint open sets
both meeting $L(G)$. Then there is a loxodromic $g$ in $G$ with a
fixed point in $O_{1}$ and a fixed point in $O_{2}$.

 \emph{\textbf{Proof.}}
First we recall that if $f$ is loxodromic with an attractive fixed
$\alpha$ and a repulsive fixed point $\beta$, then as
$n\rightarrow \infty$,$f^{n}\rightarrow \alpha$ uniformly on each
compact subgroup of $\overline{H^{n}_{C}}-\{\beta\}$ and
$f^{-n}\rightarrow \beta$ uniformly on each compact subset of
$\overline{H^{n}_{C}}-\{\alpha\}$. The repulsive fixed point of
$f$ is the attractive fixed point of $f^{-1}$.

Now consider $G$, $O_{1}$ and $O_{2}$ as in the lemma. It follows
that there is a loxodromic $p$ with attractive fixed point in
$O_{1}$ and a loxodromic $q$ with attractive fixed point in
$O_{2}$. Since $G$ is non-elementary, there is a loxodromic $f$
with attractive fixed point $\alpha$ and repulsive fixed point
$\beta$, nether fixed by $p$. Now choose and (then fix) some
sufficiently large value of $m$ so that $$g=p^{m}fp^{-m}$$ has its
attractive fixed point $\alpha_{1}(=p^{m}\alpha)$ and repulsive
fixed point $\beta_{1}(=p^{m}\beta)$ in $O_{1}$. Then choose (and
fix) some sufficiently large value of $r$ so that $$h=q^{r}$$ maps
$\alpha_{1}$ into $O_{2}$: put $\alpha_{2}=h(\alpha_{1})$.

Next, construct open convex neighborhood $E$ and $K$ of
$\beta_{1}$ and $\alpha_{2}$ with the properties
$$\beta_{1}\in E\subset \overline{E}\subset O_{1}$$
$$\alpha_{2}\in K \subset \overline{K}\subset O_{2}.$$

As $\beta_{1}$ is not in $\overline{K}$ we see that
$g^{n}\rightarrow \alpha_{1}$ uniformly on $\overline{K}$ as
$n\rightarrow \infty$. As $h^{-1}(K)$ is an open neighborhood of
$\alpha_{1}$ we see that for all sufficiently large $n$,
$$g^{n}(\overline{K})\subset h^{-1}(K)$$ and so
$$hg^{n}(\overline{K})\subset K \eqno (2.1)$$

As $h(\alpha_{1})$ is not in $ \overline{E}$ so $\alpha_{1}$ is
not in $h^{-1}(\overline{E})$ and so $g^{-n}\rightarrow \beta_{1}$
uniformly on $h^{-1}(\overline{E})$ as $n\rightarrow \infty$. Thus
for all sufficiently large $n$,
$$g^{-n}h^{-1}(\overline{E})\subset E \eqno (2.2)$$
Choose a value of $n$ for which (2.1) and (2.2) hold. By Brouwer
fixed point theorem, $hg^{n}$ is loxodromic with a fixed point in
$K $: also, $g^{-n}h^{-1}$, which is $(hg^{n})^{-1}$, has a fixed
point in $E$, hence so does $hg^{n}$. By definition, $hg^{n}$ is
not parabolic. According to Lemma 3.3.2 of [5], $hg^{n}$ is not
elliptic either. So $hg^{n}$ is a loxodromic element with one
fixed point in $O_{1}$ and the other in $O_{2}$.
\hspace{11cm}$\Box$

\textbf{Lemma 2.7.} Let $\{f_{m}\}$ be a sequence in $PU(1,n;C)$
converging to a loxodromic element $f$. Then $f_{m}$ is loxodromic
for sufficiently large $m$.

\emph{\textbf{Proof.}}  Let $x$ and $y$ be the attractive and
repulsive fixed point of $f$, respectively. We have
$$\lim_{j\rightarrow\infty} \lim_{m\rightarrow\infty}
(f_{m})^{j}(p)=\lim_{j\rightarrow\infty} f^{j}(p)=x,$$
$$\lim _{j\rightarrow\infty}\lim _{m\rightarrow\infty}(f_{m})^{-j}(p)=\lim _{j\rightarrow\infty}f^{-j}(p)=y,$$
for all $p$ in $\overline{H^{n}_{C}}\backslash \{y\}$ and
$\overline{H^{n}_{C}}\backslash \{x\}$, respectively.

Let $U,V$ be two open convex neighborhood of $x$ and $y$ in
$\overline{H^{n}_{C}}$ such that $$\overline{U}\bigcap
\overline{V}=\emptyset.$$ Then for all sufficiently large $j,m,$
$$(f_{m})^{j}(\overline{U})\subset U,\;\;\;\; (f_{m})^{-j}(\overline{V})\subset V.$$

Brouwer fixed point theorem tells us that, for all sufficiently
large $j,m$, $(f_{m})^{j}$ has one fixed point in $U$ and another
in $V$. Hence $(f_{m})^{j}$ is not parabolic. By lemma 3.3.2 in
[5], $(f_{m})^{j}$ is not elliptic either. Therefore all these
$(f_{m})^{j}$ are loxodromic. So $f_{m}$ is loxodromic for
sufficiently large $m$. In fact, for the purpose of a
contradiction, suppose $f_{m}$ is parabolic or elliptic. If
$f_{m}$ is parabolic, then $f_{m}$ has exactly one fixed point in
$\partial H^{n}_{C}$ and has $(f_{m})^{j}$, that is, $(f_{m})^{j}$
is parabolic, this is a contradiction. If $f_{m}$ is elliptic.
Then $f_{m}$ has a fixed point in $H^{n}_{C}$ and has
$(f_{m})^{j}$, that is, $(f_{m})^{j}$ is elliptic,  also a
contradiction. Consequently, $f_{m}$ must be loxodromic for
sufficiently large $m.$\hspace{7.2cm}$\Box$

We know that if $G\subset PU(1,n;C)$ is non-elementary then there
must exist infinitely many loxodromic elements in $G$. Let $h\in
G$ be some loxodromic element and let $x_{0}$ and $y_{0}$ be its
distinct fixed points. Set$$G(x_{0},y_{0})=\{f\in
G:\{x_{0},y_{0}\}\subset fix(f)\}.$$ We also need the following
lemma, which is a direct consequence of Lemma 2.2 in [12].

\textbf{Lemma 2.8.} Suppose a non-elementary subgroup $G$ of
$PU(1,n;C)$ be discrete, then $G(x_{0},y_{0})$ is a bounded
torsion.

\vskip 2pt
 \vskip 2mm {\bf 3. Discreteness Criteria for Subgroups of PU(1,$n$;C) }\vskip
 2mm

In this section, we will state our principal results. Above all,
we will introduce the first discreteness criterion for subgroups
of $PU(1,n;C)$ by using a test map which need not to be in $G$.

\textbf{Theorem 3.1.} Let $G$ be a non-elementary subgroup of
$PU(1,n;C)$ with condition $A$, and $h$ be a non-trivial element.
If each $\langle h,g\rangle$ is discrete $(g\in G)$, then $G$ is
discrete.

\emph{\textbf{Proof.}} Let $U_{i}\subset
\overline{H^{n}_{C}}(i=1,2,3)$ be disjoint open sets both meeting
$L(G)$, and $h$ does not fix any point in $U_{1}$.By lemma 2.6, we
can find loxodromic elements $f_{i}(i=1,2,3)$ in $G$ which have
the following properties:

 (i) $f_{1}$ has
its both attractive and repelling fixed points in $U_{1}$.

 (ii)
$f_{i}$ has its attractive fixed point in $U_{i}$ and repelling
fixed point in $U_{1}$ for $i=2,3.$\\
Then there is an integer $k$ such that $f_{i}^{k}(fix(h))\subset
U_{i}(i=1,2,3)$.

Suppose that $G$ satisfies the conditions of the theorem, but $G$
is not discrete. Then we can find a sequence $\{g_{j}\}$ of
distinct element in $G$ such that $g_{j}\rightarrow I$ as
$j\rightarrow \infty$. Thus we have
$$max\{N(g_{j}),N([g_{j},(f_{i}^{k}hf_{i}^{-k})^{p}]):p=1,2,...,n+1\}\rightarrow 0.$$
Since all groups $\langle
g_{j},f_{i}^{k}hf_{i}^{-k}\rangle=f_{i}^{k}\langle
f_{i}^{-k}g_{j}f_{i}^{k},h\rangle f_{i}^{-k}$ are discrete by the
assumption. In view of lemma 2.5, we get that each $\langle
g_{j},f_{i}^{k}hf_{i}^{-k}\rangle$ is elementary for large $j$.
Because $G$ satisfies Condition $A$ and $\langle
g_{j},f_{1}^{k}hf_{1}^{-k}\rangle$ is discrete, we also have
$Card(fix(g_{j}))\leq 2$ for $j\rightarrow\infty$.

 (a) $h$ is parabolic. Let $a$ be the fixed point of
$h$. We have
$$L\langle g_{j},f_{i}^{k}hf_{i}^{-k}\rangle=fix(f_{i}^{k}hf_{i}^{-k})=\{f_{i}^{k}(a)\}\;\;\;
(i=1,2,3)$$ and each $g_{j}$ fixes $f_{i}^{k}(a)\in
U_{i},\;\bigcap U_{i}=\emptyset$ $(i=1,2,3)$, this implies that
$g_{j}$ has three distinct fixed points, but
$Card(fix(g_{j}))\leq2$ for $j\rightarrow\infty$, this is a
contradiction.

(b) $h$ is loxodromic. Assume $a$ and $b$ are the fixed points of
$h$. We have$$L(\langle
g_{j},f_{i}^{k}hf_{i}^{-k}\rangle)=fix(f_{i}^{k}hf_{i}^{-k})=\{f_{i}^{k}(a),f_{i}^{k}(b)\}\;\;(i=1,2,3).$$
 $g_{j}$ either fixes both $f_{i}^{k}(a)$ and $f_{i}^{k}(b)$ or
interchanges them for sufficiently large $j$. Without loss of
generality, we may assume that for each $j$, $g_{j}$ interchanges
$f_{1}^{k}(a)$ and $f_{1}^{k}(b)$. So it follows that $g_{j}$
certainly fix both $f_{i}^{k}(a)$ and $f_{i}^{k}(b)$ $(i=2,3)$.
Since $f_{i}^{k}(a)\in U_{i}$, $f_{i}^{k}(b)\in U_{1}\;(i=2,3)$
and $U_{1}\bigcap U_{2}\bigcap U_{3}=\emptyset$, it is clear that
$g_{j}$ have at least three distinct fixed points. But
$Card(fix(g_{j}))\leq2$, as $j\rightarrow\infty$. This again leads
to a contradiction. We complete the proof of the theorem.
\hspace{8.5cm}$\Box$

\textbf{Corollary3.2.} Let $G$ be a non-elementary subgroup of
$PU(1,n;C)$, and $h\in G$ be a parabolic or loxodromic element.
Then $G$ is discrete if and only if for every element $g(\neq h)$
in $G$ the group $\langle h,g\rangle$ is discrete.\\

\textbf{Theorem 3.3.} Let a non-elementary subgroup $G$ of
$PU(1,n;C)$ satisfy condition A. Then $G$ is discrete if and only
if for two arbitrary loxodromic element $f$ and $g$ in $G$ the
group $\langle f,g\rangle$ is discrete.

\emph{\textbf{Proof.}} The necessity is obvious, we only need to
prove the sufficiency. Suppose that every two-loxodromic-generator
subgroup of $G$ is discrete and yet $G$ is not discrete. Then
there is a distinct sequence $\{g_{j}\}\subset G$ converging to
the identity. Our aim is to reach a contradiction.

As $G$ is non-elementary, there definitely exists a loxodromic
element $h$ in $G$. Since $g_{j}h\rightarrow h$ as
$j\rightarrow\infty$, it follows from Lemma 2.7 that $g_{j}h$ is
loxodromic for sufficiently large $j$. We may assume that for each
$j$, $g_{j}h$ is loxodromic. Since $h$ and $g_{j}h$ are
loxodromic, by the assumption, we know that $\langle
h,g_{j}h\rangle=\langle h,g_{j}\rangle$ is discrete. Because $G$
satisfies Condition $A$ and $\langle h,g_{j}\rangle$ is discrete,
we obtain that $Card(fix(g_{j})\leq2$ for sufficiently large $j$.

 According to Lemma 2.5 and the assumption $g_{j}\rightarrow I$ as
 $j\rightarrow\infty$, we have that $\langle h,g_{j}\rangle$ is
 discrete and elementary for sufficiently large $j$. Since $h$ is loxodromic,
 we have $g_{j}$ either fixes the fixed points of $h$ or exchanges
 them as $j\rightarrow\infty$.
  As $G$ is non-elementary, there exist another two loxodromic elements
 $f_{1},f_{2}$ such that $h\bigcap f_{1}\bigcap f_{2}=\emptyset$.
 For the above mentioned reason, it is not difficult to deduce
 that $\langle f_{i},g_{j}\rangle\; (i=1,2)$ is discrete elementary and $g_{j}$ either fixes the fixed
 points of $f_{i}\;(i=1,2)$ or interchanges them
 for enough large $j$.
Without loss of generality, we may assume that each $g_{j}$
exchanges the fixed points of $h$, so $g_{j}$ necessarily fixes
each fixed point of $f_{i}\;(i=1,2)$. However, $f_{1}$ and $f_{2}$
have no common fixed points, thus $Card(fix(g_{j}))=4$. This is a
contradiction with $Card(fix(g_{j}))\leq 2$ as
$j\rightarrow\infty$. Up to now, we complete the proof of the
theorem.\hspace{5cm}$\Box$
\\

Let $h\in G$ be some loxodromic element and let $x_{0}$ and
$y_{0}$ be its distinct fixed points. We now use $G(x_{0},y_{0})$
to strengthen theorem 3.3 as follows.

\textbf{Theorem 3.4.} Suppose that $G$ in $PU(1,n;C)$ is a
non-elementary subgroup, then $G$ is discrete if and only if

(1)$G(x_{0},y_{0})$ satisfy condition A;

 (2)every two-loxodromic-generator subgroup is discrete.

\emph{\textbf{Proof.}} In order to prove necessity, it suffices to
show that $G(x_{0},y_{0})$ has bounded torsion if $G$ is discrete.
By lemma 2.8, we know that $G(x_{0},y_{0})$ has bounded torsion.
Since a group with bounded torsion satisfies Condition A, we
directly deduce the conclusion.

Now we prove sufficiency. Suppose that $G$ is not discrete
although every subgroup generated by two loxodromic elements is
discrete. Thus there is an infinite sequence $\{g_{j}\}$ of
distinct elements such that $g_{j}\rightarrow I$ as
$j\rightarrow\infty$. We derive a contradiction as follows.

Let $h\in G(x_{0},y_{0})$ be a loxodromic element. Since
$g_{j}h\rightarrow h$ as $j\rightarrow \infty$, we get that
$\langle h,g_{j}\rangle=\langle h,g_{j}h\rangle$ is discrete for
enough large $j$, according to Lemma 2.7 and the assumption in
Theorem. As the sequence $\{g_{j}\}$ converges to the identity, we
have $$max\{N(g_{j}),N([g_{j},h^{i}]):i=1,2,...,n+1\}\rightarrow
0.$$ Thus by Lemma 2.5, there exists $J$ such that $\langle
g_{j},h \rangle$ is discrete and elementary when $j>J$. Since $h$
is a loxodromic element, $g_{j}$ fix or interchange the two fixed
points of $h$ when $j>J$, namely
$g_{j}\{x_{0},y_{0}\}=\{x_{0},y_{0}\}$. Therefore $g_{j}^{2}\in
G(x_{0},y_{0})$ as $j>J$. Since $G(x_{0},y_{0})$ satisfies
Condition $A$ and $\langle h,g_{j}^{2}\rangle$ is discrete, we
gain that $Card(fix(g_{j}^{2}))\leq2$ for sufficiently large $j$.
As $G$ is non-elementary, there exists anther two loxodromic
elements $f_{1}$ and $f_{2}$ such that $f_{i}\;(i=1,2)$ and $h$
have no common fixed points. We can also acquire that
$Card(fix(g_{j}))=4$ as $j\rightarrow\infty$, for reason see the
proof of theorem 3.3. So $g_{j}^{2}$ have at least four fixed
points as $j\rightarrow\infty$, this is a contradiction. We
complete its proof of the last theorem.\hspace{13.2cm}$\Box$

\par \hskip 20mm Department of Mathematics
\par \hskip 20mm Ocean University of China
\par \hskip 20mm Qingdao, Shandong 266071
\par \hskip 20mm P. R. China
\par \hskip 20mm E-Mail: xiaoyanbbbb@163.com
\par \hskip 20mm E-Mail: Changjunli7921@hotmail.com

\end{document}